\documentclass[10pt,twoside]{article}
\usepackage[latin1]{inputenc}
\usepackage{amsmath,amssymb}
\usepackage{Latex-document}

\renewcommand{\thesection}{\arabic{section}}
\renewcommand{\thesubsection}{\thesection.\arabic{subsection}}

\numberwithin{equation}{section}
\newcommand{\beq}{\begin{equation}}
\newcommand{\beql}[1]{\begin{equation}\label{#1}}
\newcommand{\eeq}{\end{equation}}
\newcommand{\ba}{\begin{array}}
\newcommand{\ea}{\end{array}}

\newcommand{\ts}{\textstyle}
\newcommand{\wt}{\widetilde}
\newcommand{\wh}{\widehat}
\def\dd{\;\!\mathrm{d}} 
\def\rmd{\mathrm{d}}
\def\rmw{\mathrm{w}}
\def\R{{\mathbb R}}
\def\N{{\mathbb N}}

\newcommand{\pl}{\partial}
\newcommand{\ti}{\times}
\newcommand{\reff}[1]{{\textrm{(\ref{#1})}}}
\newcommand{\set}[2]{ \{\, #1 \: | \: #2 \,\} }

\newtheorem{theorem}{Theorem}[section]
\newtheorem{definition}[theorem]{Definition}
\newtheorem{corollary}[theorem]{Corollary}

\newcommand{\PROOF}{{\bf Proof:} }
\newcommand{\QED}{\mbox{}\hfill\rule{5pt}{5pt}\medskip\par}

\def\calD{{\mathcal D}}
\def\calE{{\mathcal E}}
\def\calF{{\mathcal F}}
\def\calI{{\mathcal I}}
\def\calJ{{\mathcal J}}

\def\calR{{\mathcal R}}
\def\calS{{\mathcal S}}
\def\calZ{{\mathcal Z}}
\def\rmC{{\mathrm C}}
\def\rmD{{\mathrm D}}
\def\rmH{{\mathrm H}}
\def\rmL{{\mathrm L}}
\def\rmT{{\mathrm T}}
\def\rmW{{\mathrm W}}

\def\Diss{\mathrm{Diss}}

\newcommand{\mybox}[1]{\begin{center}
   \fbox{\begin{minipage}[c]{0.96\textwidth}\textsl{#1}%
         \end{minipage}}\end{center}}
\def\jump#1{[\![ #1]\!]_\Gamma}

\def\ext{{\mathrm{ext}}}
\def\red{{\mathrm{red}}}
\def\old{{\mathrm{old}}}
\def\new{{\mathrm{new}}}
\title{\bf Analysis of Energetic Models for\vskip -2mm Rate-Independent Materials\vskip 6mm}
\author{Alexander Mielke\vspace*{-0.5cm}\thanks{Mathematisches Institut A,
Universität Stuttgart, Pfaffenwaldring 57, 70569 Stuttgart,
Germany. E-mail: mielke@mathematik.uni-stuttgart.de}}
\date{\vspace{-8mm}}

\begin{document}

\maketitle

\thispagestyle{first} \setcounter{page}{817}

\begin{abstract}

\vskip 3mm

We consider rate-independent models which are defined via two
functionals: the time-dependent energy-storage functional
$\calI:[0,T]\ti X\to [0,\infty]$ and the dissipation distance
$\calD:X\ti X\to[0,\infty]$. A function $z:[0,T]\to X$ is called a
solution of the {energetic model}, if for all $0\leq s<t\leq T$ we
have
\\[1mm]
\begin{tabular}{@{}r@{\;\;}c}{stability:}&
    $\mathcal I(t,z(t)) \leq \mathcal I(t,\widetilde z)+\calD(z(t),\wt z)$
    for all $\wt z\in X$;\\
    {energy inequality:}& $\mathcal I(t,z(t)) {+} \Diss_\calD(z,[s,t]) \leq
    \mathcal I(s,z(s)) {+} \int_s^t \partial_\tau\mathcal
    I(\tau,z(\tau))\,\mathrm d \tau$.
\end{tabular}
\smallskip

We provide an abstract framework for finding solutions of this
problem. It involves time discretization where each incremental
problem is a global minimization problem. We give applications in
material modeling where $z\in \calZ\subset X$ denotes the internal
state of a body. The first application treats shape-memory alloys
where $z$ indicates the different crystallographic phases. The
second application describes the delamination of bodies glued
together where $z$ is the proportion of still active glue along
the contact zones.
 The third application treats
finite-strain plasticity where $z(t,x)$ lies in a Lie group.

\vskip 4.5mm


\end{abstract}

\vskip 12mm

\section{Introduction} \label{s1}

\vskip-5mm \hspace{5mm}

Many evolution equations can be written in the abstract form
\beql{e1.1}
 0\in \pl \Psi( \dot z(t)) +\mathrm D \calI(t,z(t)),
\eeq where $z\in X$ is the state variable, $\calI$ is the
energy-storage functional, $\Psi:X\to[0,\infty]$ is a convex
dissipation functional, and  $\pl\Psi$ means the set-valued
subdifferential (see \cite{ColVis90CDNE} for this doubly nonlinear
form). Rate-independency is realized by assuming that $\Psi$ is
homogeneous of degree 1.

We replace the above differential inclusion by a weaker energetic
formulation, which is also more general since it allows for
$z$-dependent dissipation functionals. For given $\calI:[0,T]\ti
X\to[0,\infty]$ and a given dissipation distance $\calD:X\times X
\to [0,\infty]$ satisfying the triangle inequality, we impose the
energetic conditions of \emph{global stability} (S) and the
\emph{energy inequality} (E) instead of \reff{e1.1}. A function
$z:[0,T]\to X$ is called a \emph{solution of the energetic model},
if for all $0\leq s<t\leq T$ we have
\\[1mm]
\hspace*{10mm} (S) \quad
$\calI(t,z(t))\leq \calI(t,\wt z)+\calD( z(t),\wt z)$ for all $\wt z\in X$;
\\
\hspace*{10mm} (E) \quad
$\calI(t,z(t))+\Diss_\calD(z,[s,t]) \leq \calI(s,z(s))+\int_s^t
\pl_\tau\calI(\tau,z(\tau))\dd \tau$.
\\[1mm]
Here,  $\Diss_\calD(z,[s,t])$ is called the dissipation of $z$ on the interval
$[s,t]$ and is defined as the supremum of
$\sum_{j=1}^N \calD(z(t_{j-1}),z(t_j))$ over all $N\in \N$ and all
discretizations $s=t_0<t_1<\ldots <t_N=t$.

Assuming $\calD(z_0,z_1)=\Psi(z_1{-}z_0)$, convexity of
$\calI(t,\cdot)$ and further technical assumptions, this energetic
formulation is equivalent to \reff{e1.1}, see
\cite{MieThe01?RIHM}. However, the latter form is more general as
it applies to nonconvex problems and it doesn't need
differentiability of $t\mapsto z(t)$ nor of $z\mapsto \calI(t,z)$.
A related energetic approach to equations of the type \reff{e1.1}
is presented in \cite{Visi01NAE}, however, it remains unclear
whether that method applies to the rate-independent case.

In Section \ref{s2} we discuss the abstract setting in more detail
and in Section \ref{s:3} we provide existence results for
solutions for given initial values $z(0)=z_0$. The existence
theory is based on time-incremental minimization problems of the
form
\[
z_k \in \mathrm{argmin}\set{ \calI(t_k,z){+}\calD(z_{k-1},z) }{ z \in X}
\]
and the BV bound for $z:[0,T]\to X$ obtained via the dissipation functional
satisfying $\calD(z_0,z_1) \geq c_D \left\| z_0{-}z_1\right\|$.
However, one needs additional compactness properties, if $X$ is infinite
dimensional. Here we propose a version where $\calI$ satisfies coercivity with
respect to an embedded Banach space $Y$, i.e., $\calI(t,z)\geq -C_1{+}c_1
\|z\|_Y^\alpha $ with $c_1,C_1,\alpha>0$, where $Y$ is compactly embedded in
$X$.

For the case of $\calD$ having the form $\calD(z_0,z_1)=\Psi(z_1{-}z_0)$ this
theory was developed in \cite{MieThe01?RIHM}.  The case of general $\calD$
can be found in \cite{MaiMie02?}.
\medskip

The flexibility of the energetic formulation allows for applications in
continuum mechanics, where $z: \Omega \to Z$ plays the r\^ole of internal
variables in the material occupying the body $\Omega \subset \R^d$. Note that
$Z$ may be a manifold containing the internal variables like phase indicators,
plastic or phase transformations, damage, polarization or magnetization. By
$\calZ$ we denote the set of all admissible internal states.  The elastic
deformation is $\varphi : \Omega \to \R^d$ and $\calF$ denotes the set of
admissible deformations $\varphi$.

Energy storage is characterized via the functional $\calE: [0,T]
\ti \calF \ti \calZ \to \R$, where $t\in [0,T]$ is the
(quasi-static) process time, which drives the system via changing
loads. In typical material models, $\calE$ has the form
\[\ts
\calE(t,\varphi,z)=\int_\Omega W(x,\rmD\varphi(x),z(x))\dd x -\langle
\ell_\ext (t),\varphi\rangle,
\]
where $W$ is the stored-energy density and $\ell_\ext(t)$ denotes
the external loadings.

Dissipation is characterized by an infinitesimal
Finsler metric $\Delta : \Omega \ti \rmT Z \to [0,\infty]$,
such that the curve $z:[t_0, t_1] \to \calZ$ dissipates the energy
\[\ts
\mathrm{Diss}\, (z, [t_0, t_1]) =
\int^{t_1}_{t_0}\int_\Omega \Delta (x, z(t, x), \dot{z} (t, x))\dd x \dd t.
\]
The global dissipation distance $\calD(z_0,z_1)$ is then the infimum
over all curves connecting $z_0$ with $z_1$.
The relation to the abstract theory above is obtained by eliminating
the elastic deformation  via
\[\ts
\calI(t,z)=\inf\{\: \calE(t,\varphi,z)\:|\: \varphi\in \calF\:\}
\mbox{ for } z \in \calZ \mbox{ and } \calI (t,z) = + \infty \mbox{ else}.
\]
Obviously, the functional $\calI$ is now fairly complicated and it is
important to have rather general conditions in the abstract theory.
\medskip

In Section \ref{s:cont} we illustrate the usefulness of the
abstract approach by discussing three quite different
applications; however, the theory is used in other areas as well,
e.g., in fracture mechanics \cite{FraMar98RBFE,DalToa02MQSG} and
in micro-magnetics \cite{Kruz02VMMS,RouKru02?MEMM}.

Our first model describes phase transformations in shape-memory
alloys as discussed in
\cite{MieThe99MMRI,MiThLe02VFRI,Roub02?EMMP,GoMiHa02FEMV}.  Here
$z:\Omega\to Z$ indicates either the microscopic distribution of
the phases or a mesoscopic average of the microscopic
distribution. In the first case we choose
$Z=Z_p=\{e_1,\ldots,e_p\}\subset\R^p$, where $e_j$ denotes the
$j$-th unit vector in $\R^p$ and in the second case we choose
$Z=\mathrm{conv}\,Z_p$. In both cases the dissipation distance is
given by a volume integral measuring the amount of volume which is
transformed into another phase: $\calD(z_0,z_1)=\int_\Omega
\Delta(z_1(x){-}z_0(x))\dd x$, where $\Delta:\R^p\to
\left[0,\infty\right[$ is convex and homogeneous of degree $1$.
This leads naturally to the basic space $X=\rmL^1(\Omega,\R^p)$
and $\calZ=\set{ z\in X}{ z(x)\in Z\text{ a.e.}}$.

Including in $\calE$ an interfacial energy proportional to the area of
the interfaces between regions of different phases provides a reduced energy
$\calI$ which is coercive in $Y=\mathrm{BV}(\Omega,\R^p)$, see \cite{Main02?}.
For an existence result in the  case without interfacial energy we refer to
\cite{MiThLe02VFRI}.
\smallskip

The second application describes the delamination of a body $\Omega$ which
is glued together along $n$ hypersurfaces
$\Gamma_j$, $j=1,\ldots,n$.   The internal state $z:\Gamma=\cup_1^n
\Gamma_j\to [0,1]$ denotes the percentage of glue along $\Gamma$ which
remains in effect. The dissipation is given by a material constant $c_\calD$
times the destroyed glue, i.e.,  $\calD(z_0,z_1)= c_\calD \int_\Gamma
z_0(x){-}z_1(x) \dd a(x)$ for $z_1 \leq z_0$ and $\calD (z_0,z_1)=+
\infty$ else. The basic underlying space is
$\rmL^1(\Gamma)$ and now compactness arises via the trace operator
$\rmH^1(\Omega) \to \rmL^2(\Gamma)$ which makes the reduced energy functional
$\calI$ weakly continuous.
\smallskip

The final application is devoted to the modeling of
elasto-plasticity with finite strains. There the internal variable
$z=(P,p)$ consists of the plastic transformation $P\in
\mathrm{SL}(d)$ and hardening parameters $p\in \R^k$. Invariance
under previous plastic deformations leads to dissipation metrics
which are left-invariant, i.e., $\Delta((P,p),(\dot P,\dot p))=
\Delta((I,p),(P^{-1}\dot P,\dot p))$. This geometric nonlinearity
clearly shows that we need general dissipation distances $\calD$
avoiding any linear structure. In single-crystal plasticity
$\Delta$ is piecewise linear in $P^{-1}\dot P\in \mathrm{sl}(d)$
which leads to Banach manifolds and the dissipation metric is then
a left-invariant Finsler metric. For applications in this context
see \cite{CaHaMi02NCPM,Miel01?FELG,Miel02?EFME}.

\section{Abstract setup of the problem}
\label{s2}

\vskip-5mm \hspace{5mm}

We start with a Banach space $X$ which is not assumed to be reflexive, since
our applications in continuum mechanics (cf.\ Section \ref{s:cont})
naturally lead to spaces of the form $\rmL^1(\Omega,\R^k)$.
The first ingredient of the energetic formulation
is the \emph{dissipation distance}
$\calD:X\ti X\to [0,\infty]$ satisfying the triangle inequality:
\[
\calD(z_1,z_3) \leq \calD(z_1,z_2) + \calD(z_2,z_3)\quad
\text{for all }z_1,z_2,z_3\in X.
\]
We don't enforce symmetry, i.e., we allow for $\calD(z_0,z_1) \neq
\calD(z_1,z_0)$ as in Section \ref{ss:delam}. We assume that there
is a constant $c_\calD>0$ such that $\calD(z_0,z_1)\geq
c_\calD\|z_1{-}z_0\|_X$ for all $z_0, z_1\in X$. The latter
condition is in fact the one which determines the appropriate
function space $X$ for a specific application. Moreover, $\calD$
is assumed to be s-weakly lower semicontinuous. (We continue to
use the abbreviation s-weak for ``sequentially weak''.) We call
$\calD(z_0,z_1)$ the dissipation distance from $z_0$ to $z_1$.

For a given curve $z:[0,T]\to X$ we define the total dissipation on $[s,t]$
via
\beql{e:Diss}\ts
\Diss_\calD(z;[s,t])=\sup\set{
   \sum_1^N \calD(z(\tau_{j-1}),z(\tau_j)) }{ N {\in} \N,
s{=}\tau_0{<}\tau_1{<}\cdots{<}\tau_N{=}t}.
\eeq

The second ingredient is the energy-storage functional
$\calI:[0,T]\ti X\to [0,\infty]$, which is assumed to be bounded
from below and then normalized such that it takes only nonnegative
values. Here $t\in [0,T]$ plays the r\^ole of a (very slow)
process time which changes the underlying system via changing
loading conditions. For fixed time $t$, the map
$\calI(t,\cdot):X\to [0,\infty]$ is assumed to be s-weakly lower
semicontinuous, i.e., $z_j\rightharpoonup z$ implies
$\calI(t,z)\leq \liminf_{j\to \infty} \calI(t,z_j)$. Moreover, we
assume that for all $z$ with $\calI(t,z)<\infty$ the function
$t\mapsto \calI(t,z)$ is Lipschitz continuous with $| \pl_t
\calI(t,z)| \leq C_\calI$.

\begin{definition}
A curve $z:[0,T]\to X$ is called a \textbf{solution} of the
rate-independent model $(\calD,\calI)$, if \textbf{global
stability (S)} and \textbf{energy
  inequality (E)} holds:
\mybox{%
\textbf{(S)} For all $t\in [0,T]$ and all $\wh z\in X$ we have $\calI(t,z(t))
\leq I(t,\wh z)+\calD(z(t),\wh z)$.
\\[0.3em]
\textbf{(E)} For all $t_0,t_1$ with $0\leq t_0<t_1\leq T$ we have\\[0.2em]
\centerline{$
\calI(t_1,z(t_1))+\Diss_\calD
(z;[t_0,t_1])\leq \calI(t_0,z(t_0))+\int_{t_0}^{t_1}
\pl_t \calI(t,z(t))\dd t.
$}
}
\end{definition}

The definition of solutions of (S)\&(E) is such that it implies the two
natural requirements for evolutionary problems, namely that
\emph{restrictions}  and \emph{concatenations} of solutions remain
solutions. To be more precise, for any solution $z:[0,T] \to E$ and
any subinterval $[s,t] \subset [0,T]$, the restriction $z|_{[s,t]}$
solves (S)\&(E) with initial datum $z(s)$. Moreover, if $z_1:[0,t] \to E$
and $z_2:[t,T] \to E$ solve (S)\&(E) on the respective intervals and if
$z_1(t) = z_2(t)$, then the concatenation $z:[0,T] \to E$ solves (S)\&(E)
as well. Under a few additional assumptions, it is shown in
\cite{MieThe01?RIHM} that (S) and (E) together imply
that, in fact,  the energy inequality is in an equality, i.e., for $0\leq t_0
<t_1\leq T$ we have
\begin{equation}
 \label{eq:2.5}\ts
  \calI(t_1,z(t_1)) + \Diss_\calD(z;[t_0,t_1]) =
  \calI(t_0,z(t_0)) + \int_{t_0}^{t_1} \pl_t \calI(t,z(t)) \dd t.
\end{equation}

Rate-independency manifests itself by the fact that the problem
has no intrinsic time scale. It is easy to show that  $z$ is a
solution for $(\calD,\calI)$ if and only if the reparametrized
curve $\wt z:t\mapsto z(\alpha(t))$, with $\dot\alpha > 0$, is a
solution for $(\calD,\wt \calI)$, where
$\wt\calI(t,z)=\calI(\alpha(t),z)$. In particular, the stability
(S) is a static concept and the energy estimate (E) is
rate-independent, since the dissipation defined via \reff{e:Diss}
is scale invariant like the length of a curve.

The major importance of the energetic formulation is that neither the given
functionals $\calD$ and $\calI(t,\cdot)$ nor the solutions $z:[0,T] \to X$
need to be differentiable. In particular, applications in continuum
mechanics often have low smoothness. Of course, under additional
smoothness assumptions on $\calD$ and $\calI$ the weak energetic form (S)\&(E)
can be replaced by local formulations in the form of differential inclusions
like \reff{e1.1} (\cite{ColVis90CDNE,Visi01NAE})
or variational inequalities. See \cite{MieThe01?RIHM} for a discussion of the
implications between these different formulations.


\section{Time discretization and existence} \label{s:3}

\vskip-5mm \hspace{5mm}

The major task is now to develop an existence theory for the initial value
problem, i.e., to find a solution in the above sense which additionally
satisfies $z(0)=z_0$. In general, we should not expect uniqueness without
imposing further conditions like smoothness and uniform convexity of
$\calI(t,\cdot)$ and $\calD$, see \cite{MieThe01?RIHM}.

The stability condition (S) can be rephrased by defining the stable sets
\[\ts
\calS(t) := \set{z \in X }{ \calI(t, z) \leq \calI(t, \wh z) + \calD (z,\wh z)
\text{  for all } \wh z \in X }.
\]
Then, (S) simply means $z(t)\in \calS(t)$ for all $t\in [0,T]$. The properties
of the stable sets turn out to be crucial for deriving existence results.

One of the standard methods to obtain solutions of nonlinear evolution
equations is that of approximation by time discretizations. To this
end we choose discrete times $0 = t_0 < t_1 < \ldots < t_N = T$ and seek
$z_k$ which approximates the solution $z$ at $t_k$, i.e., $z_k\approx z(t_k)$.
Our energetic approach has the major advantage that the values $z_k$ can be
found incrementally via minimization problems. Since the methods of
the calculus of variations are especially suited for applications in material
modeling this will allow for a rich field of applications.

To motivate the following incremental variational problem consider the
nonlinear parabolic problem $h(\pl_t u)=\mathrm{div}(A\;\!\rmD u)+g $,
where we assume $h'(v)\geq 0$. The associated
fully implicit incremental problem  reads
\[\ts
h(\frac{1}{t_k{-}t_{k-1}}(u_k{-}u_{k-1}))=\mathrm{div}(A\;\!\rmD
u_k)+g(t_k).
\]
With $H(v)=\int_0^vh(w)\dd w$ we see that $u_k$ must be a
minimizer of the functional
\[\ts
\calJ_k(u_{k-1};\cdot): u\mapsto \int_\Omega (t_k{-}t_{k-1})H
(\frac{1}{t_k-t_{k-1}}
(u{-}u_{k-1})) +\frac12 \langle A\;\!\rmD u,\rmD u\rangle -g(t_k)u \dd x.
\]
In the simplest rate-independent case the function $h$ is given by
the signum function which implies $H(v)=|v|$. Hence, the length
$t_k{-}t_{k-1}$ of the $k$-th time step disappears in the
functional $\calJ_k$.  In our more general setting the incremental
problem takes the following form:

\mybox{{\bf (IP)} \sl
For $z_0 \in X$ with $\calI(0,z_0)<\infty$  find $z_1,\ldots,z_N \in X$
such that
\begin{equation}\label{eq:3.1}
  z_k \in \mathop{\mathrm{argmin}}
 \set{ \calI(t_k,z) + \calD(z_{k-1},z)}{z\in X}
\quad\text{ for }k = 1,\ldots,N.
\end{equation}
}

Here ``argmin'' denotes the set of all minimizers. Using the
s-weak lower semi-continuity of $\calD$ and $\calI$ and the
coercivity $\calI(t,z){+}\calD(z_{k-1},z)\geq
c_\calD\|z{-}z_{k-1}\|$ we obtain the following result.

\begin{theorem} \label{t:3.1}
The incremental problem \reff{eq:3.1} always has a solution. Each
solution satisfies, for $k = 1,\ldots,N$, the following properties:
\\[0.5ex]
(i)\quad $z_k$ is stable for time $t_k$, i.e., $z_k\in \calS(t_k)$;
\\[0.5ex]
(ii)\quad  $ \int_{[t_{k-1},t_k]}\pl_s \calI(s,z_k) \dd s
      \leq \calI(t_k,z_k) - \calI(t_{k-1},z_{k-1}) + \calD(z_{k-1},z_k)$ \\
\phantom{(ii)\quad  $ \int_{[t_{k-1},t_k]}\pl_s \calI(s,z_k) \dd s$}
  $\leq \int_{[t_{k-1},t_k]} \pl_s \calI (s,z_{k-1}) \dd s$;
\\[0.5ex]
(iii)\quad $\calI(t_k,z_k)+\sum^k_{j=1} \calD (z_{j-1}, z_j) \leq \calI(0,z_0)
+ C_\calI T$;  \\[0.5ex]
(iv)\quad $\|z_k\| \leq \|z_0\| {+}(\calI(0,z_0) {+} C_\calI T)/ c_\calD $.
\end{theorem}
The assertions (i) and (ii) are the best replacements for the
conditions (S) and (E) in the time-continuous case.
\medskip

For each discretization $P=\{0,t_1,\ldots,t_{N-1},T\}$ of the
interval $[0,T]$ and each incremental solution
$(z_k)_{k=1,\ldots,N}$ of (IP) we define two piecewise constant
functions which attain the values $z_k$ at $t_k$ and are constant
in-between: $Z^P$  is continuous from the left and $\wh Z^P$ is
continuous from the right. Summing the estimates (ii) in Theorem
\ref{t:3.1} over $k=j,\ldots, m$ we find the following two-sided
energy estimate.

\begin{corollary} \label{c:3.2}
Let $P$ be any discretization of $[0,T]$ and $(z_k)_{k=0,\ldots,
N}$ a solution of (IP), then for $0\leq j<m\leq N$ we have the
two-sided energy inequality
\[
\ba{r@{\:}c@{\:}l}
\calI(t_j,Z^P(t_j))+ \int_{t_j}^{t_m}\pl_s \calI(s,Z^P(s))\dd s
 &\leq& \calI(t_m,Z_P(t_m))+\Diss_\calD (Z_P,[t_j,t_m]) \\
&\leq& \calI(t_j,Z^P(t_j))+\int_{t_j}^{t_m} \pl_s \calI(s,\wh Z_P(s))\dd s.
\ea
\]
\end{corollary}

The existence of solutions can now be established by taking a sequence
$(P{(l)})_{l\in \N}$ of discretizations whose fineness
$\delta^{(l)}=\max\set{t_j^{(l)}{-}t_{j-1}^{(l)} }{j=1,\ldots,N^{(l)}}$ tends
to $0$. Moreover we assume that the sequence is hierarchical with
$P{(l)}\subset P{(l{+}1)}$. The associated solutions of (IP)$^{(l)}$ define
$z^{(l)}:=Z^{P{(l)}}$. The construction of a solution of (S)\&(E) consists now
of two parts.

First we use the dissipation bound (iii) of Theorem
\ref{t:3.1} to obtain an a priori bound in BV$([0,T],X)$:
\[\ts
c_\calD \int_{[0,T]} \|\rmd z^{(l)}\|_X \leq \Diss_\calD(z^{(l)},[0,T]) \leq
\calI(0,z_0)+C_\calI T.
\]
Then, Part (iv) in Theorem \ref{t:3.1} and the following additional compactness condition
\reff{eq:3sc} allows us to apply a selection principle in the spirit of Helly.
\beql{eq:3sc}
\ba{l}
\text {For all $R>0$ and all $t\in [0,T]$ the sets} \\
\calR_R^t:=\set{z\in X}{
  \calD(z_0,z)\leq R,\ \calI(t,z)\leq R }\text{ are s-weakly compact.}
\ea
\eeq
Thus, we can extract a
subsequence $(l_n)_{n\in \N}$ such that for all $t\in [0,T]$ the sequence
$z^{(l_n)} (t),n\in \N$, converges weakly to a limit $z^{(\infty)}(t)$ with
$\Diss_\calD(z^{(\infty)},[0,T])\leq \liminf_{n\to \infty}
\Diss_\calD(z^{(l_n)},[0,T])$.

Second we need to show that $z^{(\infty)}$ is a solution of (S)\&(E).
Using Corollary \ref{c:3.2} it is easy to give conditions which guarantee that
$z^{(\infty)}$ satisfies (E) for $t_0=0$ and $t_1=T$, and by \reff{eq:2.5}
this is sufficient. To obtain stability of $z^{(\infty)}$ there are
essentially two different ways. If additional compactness properties
allow us to
conclude that the convergence of $z^{(l_n)}(t)$ to $z^{(\infty)}$ also happens
in the strong topology, then we are in the good case. Then it suffices to know
that the set
\[
\calS_{[0,T]}=\set{(t,z)\in [0,T]\ti X}{z\in \calS(t)}=\cup_{t\in[0,T]}
(t,\calS(t))
\]
is closed in the strong topology. If strong convergence cannot be
deduced, one needs to show that $\calS_{[0,T]}$ is s-weakly
closed. This property is quite hard to obtain, since even under
nice convexity assumptions on $\calI(t,\cdot)$ the sets $\calS(t)$
are generally not convex.

The following theorem provides two alternative sets of assumptions which
enables us to turn the above construction into a rigorous existence proof.

\begin{theorem}\label{t:3.3}
  Let $\calD$ and $\calI$ be given as above and satisfy \reff{eq:3sc}.  If
  one of the conditions (a) or (b) is satisfied, then for each $z_0\in X$ with
  $\calI(0,z_0)<\infty$ there is at least one solution $z\in
  \mathrm{BV}([0,T],X)$ of (S) $\&$ (E) with $z(0)=z_0$.
  \\[0.3em]
  (a) The set $\calS_{[0,T]}$ is s-weakly closed and
  $z\mapsto \pl_t \calI(t,z)$ is s-weakly continuous.
  \\[0.3em]
  (b) The sets $\calR_R^t$ in \reff{eq:3sc} are compact, the set
  $\calS_{[0,T]}$ is closed, and $z\mapsto \pl_t \calI(t,z)$ is continuous
  (all in the norm topology of $X$).
\end{theorem}

Simple nontrivial applications of this theorem with either condition (a) or
(b) are as follows: Let $X=\rmL^1(\Omega)$ with  $\Omega\subset \R^d$
bounded and choose the dissipation distance
$\calD(z_0,z_1)=c_\calD\|z_1{-}z_0\|_X=c_\calD
\int_\Omega|z_1(x){-}z_0(x)|\dd x$. As a first case consider
\[\ts
\calI_1(t,z)=\int_\Omega \alpha(x)|z(x)|^\beta{-}g(t,x)z(x)\dd x+\gamma ,
\]
where $\alpha(x)\geq \alpha_0>0$, $\beta > 1$, and $g\in
\rmC^1([0,T],\rmL^\infty(\Omega))$. The sets $\calR_R^t$ are
closed convex sets which lie in the intersection of an
$\rmL^1$-ball and an $\rmL^\beta$-ball. Hence, we obtain the
s-weak compactness condition \reff{eq:3sc}.  Yet, $\calR^t_R$ is
not strongly compact in $\rmL^1(\Omega)$. The stable sets for
$\calI_1$ are given by
\[\ts
\calS_1(t)= \set{z\in \rmL^1(\Omega)}{
|z(x)|^{\beta-2}z(x) \in [\frac{g(t,x)-c_\calD}{\alpha(x)\beta},
\frac{g(t,x)+c_\calD}{\alpha(x)\beta}] \text{ for a.a.\:}x\in \Omega },
\]
which shows that they are s-weakly closed since they are convex
and closed. Hence, condition (a) is satisfied.

As a second case consider the nonconvex energy functional
\[\ts
\calI_2(t,z)=\int_\Omega \frac12 |\rmD z(x)|^2 {+}f(t,x,z(x))\dd x \text{ for
  }  z\in \rmH^1(\Omega)\quad\text{and }
+\infty \text{ else},
\]
where $f:[0,T]\ti \Omega \ti \R\to \R$  and $\pl_t f$ are
continuous and bounded. Now, $\calR_R^t$ is already compact in
$\rmL^1(\Omega)$ since it is closed and contained in an
$\rmH^1$-ball. With these properties, it can be shown that
condition (b) of Theorem \ref{t:3.3} holds.

\section{Applications in continuum mechanics}
\label{s:cont}

\vskip-5mm \hspace{5mm}

The flexibility of the energetic formulation allows for applications in
continuum mechanics. We consider an elastic body which is given through a
bounded domain $\Omega \subset \R^d$ with sufficiently smooth boundary.
The elastic deformation is given
by the mapping $\varphi : \Omega \to \R^d$, and the set of all admissible
deformations is denoted by $\calF$, which implements the displacement boundary
conditions.

The variable $z\in Z$ includes all the internal variables like
phase indicators, plastic or phase transformations, damage,
polarization or magnetization. A function $z: \Omega \to Z$ gives
the internal state of the material, and $\calZ$ denotes the set of
all admissible internal states. Note that $Z$ may be a manifold
with (nonsmooth) boundary. In plasticity we have
$Z=\mathrm{SL}(d)\ti \R^k$, in phase transformations we let
$Z=\set{ z\in [0,1]^k}{ \sum_1^p z^{(j)}=1}$, and in
micro-magnetism $z$ is the magnetization satisfying $|z(t,x)| =
m_0 > 0$. Moreover, below we will also consider an application
where $z$ is not defined on all of $\Omega$ but at certain parts
of the boundary.

Energy storage is characterized via the functional $\calE: [0,T]
\ti \calF \ti \calZ \to \R$ which is the sum of the total elastic energy and
the potential energies due to exterior loadings (Gibbs' energy):
\[\ts
\calE(t,\phi,z)=\int_\Omega W(x,\rmD \phi(x),z(x))\dd x -\langle
\ell_\ext(t),\phi \rangle .
\]
Here $t\in [0,T]$ is the (quasi-static) process time which drives
the system and the external loads are $\langle\ell_\ext(t),\phi
\rangle=\int_\Omega f_\ext(t,x)\cdot \phi(x)\dd x
+\int_{\Gamma_{\mathrm{tract}}} g_\ext(t,x)\cdot \phi(x)\dd a(x)$.

Dissipation is characterized via the metric $\Delta : \Omega \ti \rmT Z \to [0,
\infty]$ such that the curve $z:[0,T] \to \calZ$ dissipates the energy
\[\ts
\mathrm{Diss}\, (z, [t_0, t_1]) =
\int^{t_1}_{t_0}\int_\Omega \Delta (x, z(t, x), \dot{z} (t, x))\dd x\dd t
\quad\text{ on }[t_1, t_2].
\]
For each material point $x\in \Omega$, the infinitesimal metric
$\Delta(x,\cdot,\cdot):\rmT Z\to [0,\infty]$
defines a global distance function
$D(x,\cdot,\cdot):Z\ti Z\to [0,\infty]$ and on $\calZ$ we obtain
the global dissipation distance
\[\ts
\ba{l}
\calD(z_0,z_1)=\int_\Omega D(x,z_0(x),z_1(x))\dd x\\
=\inf\set{\mathrm{Diss}\, (z, [t_0, t_1])}{ z\in \rmC^\mathrm{Lip}([0,1],
\calZ),  z(0)=z_0,z(1)=z_1} .
\ea
\]

The rate-independent problem for this material model is defined as
in the above abstract part, but now the elastic deformation
appears as an additional variable, which, however, does not
generate any dissipation.

\begin{definition} \label{d:4.1}
A pair $(\phi,z):[0,T]\to \calF\ti \calZ$ is called a
\textbf{solution} of the rate-independent problem associated with
$\calD$ and $\calE$ if the \textbf{global stability (S)} and the
\textbf{energy inequality (E)} hold:
\mybox{%
\textbf{(S)} For all $t\in [0,T]$ and all $(\wh\phi,\wh z)\in \calF\ti \calZ$
we have\\[0.2em]
\centerline{$
\calE(t,\phi(t),z(t)) \leq \calE(t,\wh\phi,\wh z)+\calD(z(t),\wh z).
$}\\[0.3em]
\textbf{(E)} For all $t_0,t_1$ with $0\leq t_0<t_1\leq T$ we have\\[0.2em]
\centerline{$
\calE(t_1,\phi(t_1),z(t_1)){+}\Diss_\calD
(z;[t_0,t_1])\leq \calE(t_0,\phi(t_0),z(t_0)){+}\int_{t_0}^{t_1}
\pl_t \calE(t,\phi(t),z(t))\dd t.
$}
}
\end{definition}

The connection with the above abstract theory
is obtained by minimization with respect
to the deformations $\phi\in \calF$, since the
stability condition implies that $\phi(t)$ must be a minimizer of
$\calE(t,\cdot ,z(t))$. We define the associated
$\calI$  via
\[\ts
\calI(t,z)=\inf\set{ \calE(t,\varphi,z)}{ \varphi\in \calF }\text{ for } z\in
\calZ \quad \text{and }+\infty\text{ else}.
\]
While this elimination is suitable for an abstract treatment, the practical
approximation of solutions via the incremental approach is better done by
keeping the deformation and eliminating the internal variable in each
incremental step. In fact, in (IP) we now have to find
\beql{eq:4.IP}
(\phi_k,z_k)\in \mathop{\mathrm{argmin}} \set{\calE(t_k,\wh\phi,\wh
  z)+\calD(z_{k-1},\wh z) }{ (\wh\phi,\wh z)\in \calF\ti \calZ}.
\eeq
In this minimization problem the internal variable occurs only locally under
the integral over $\Omega$ and hence can be eliminated pointwise.
Defining the local reduced constitutive functions
\beql{e4.2}
\ba{l}\Psi^\red(z_\old;x,F):=\min\set{ W(x,F,z)+D(x,z_\old,z)}{z\in Z}, \\
Z_\new (z_\old;x,F)\in \mathop{\mathrm{argmin}} \set{ W(x,F,z)+D(x,z_\old,z)}{z\in Z},
\ea
\eeq
and the reduced functional $\calE^\red(z_\old;t,\phi)=\int_\Omega \Psi^\red
(z_\old;\rmD \phi)\dd x -\langle \ell_\ext(t),\phi\rangle$
the solution of \reff{eq:4.IP} is equivalent to finding $\phi\in
\mathrm{argmin} \set{\calE^\red(z_{k-1};t_k,\wh \phi)}{\wh \phi \in \calF}$
and then letting $z_k=Z_\new (z_{k-1};\rmD\phi_k)$. For more details we
refer to \cite{Miel02?EFME}.

\subsection{Phase transformations in shape-memory alloys}
\vskip -5mm \hspace{5mm}

We assume that, in each microscopic point $y$, an elastic material
is free to choose one of $p$ crystallographic phases and that the
elastic energy density $W$ is then given by $W_j(\rmD \phi)$. If
the model is made on the mesoscopic level, then the internal
variables are phase portions $z^{(j)}\in [0,1]$ for the $j$-th
phase. We set $Z=\set{z\in [0,1]^p\subset
  \R^p}{\sum_1^pz^{(j)}=1}$ and $X=\rmL^1(\Omega,\R^p)$. The material
properties are described by a mixture function $W:\R^{d\ti d}\ti
Z\to[0,\infty]$, see \cite{Miel00EMFM,MiThLe02VFRI,GoMiHa02FEMV}. The
dissipation can be shown to have the form $D(z_0,z_1)=\psi(z_1{-}z_0)$ with
$\psi(v)= \max\set{ \sigma_m\cdot v}{m=1,\ldots,M}\geq C_\psi|v|$, where
$\sigma_m\in \R^p$ are thermodynamically conjugated threshold values.

So far we are unable to prove existence results for this model in
its full generality. However, the case with only two phases
$(p=2)$ has been treated in \cite{MiThLe02VFRI} under the
additional assumption that the elastic behavior is linear and both
phases have the same elastic tensor. In that case, one sets
$z=(\theta,1{-}\theta)$ with $\theta\in [0,1]$. It can be shown
that $\calI$ is a quadratic functional in $\theta \in
\rmL^1(\Omega,[0,1])\subset \rmL^2(\Omega)$. It then follows that
the compactness condition \reff{eq:3sc} holds and condition (a) in
Theorem \ref{t:3.3} can be verified using the H-measure to handle
the weak convergence of the nonconvex terms.
\medskip

A microscopic model is treated in \cite{Main02?}. There no phase
mixtures are allowed, i.e., we assume $z\in
Z_p:=\{e_1,e_2,\ldots,e_p\}\subset \R^p$, where $e_j$ is the
$j$-th unit vector. Thus, the functions $z\in \calZ$ are like
characteristic functions which indicate exactly one phase at each
material point. The dissipation is assumed as above, but now the
elastic energy contains an additional term measuring the surface
area of the interfaces between the different regions:
\[\ts
\calE(t,\phi,z)=\int_\Omega W(\rmD\phi,z)\dd x +\sigma \int_\Omega |\rmD z|
-\langle \ell_\ext(t),\phi\rangle ,
\]
where $\sigma$ is a positive constant and $\int_\Omega |\rmD z|$ is
$\sqrt2$ times the area of all interfaces. Here $\calZ=\set{z: \Omega\to Z_p}{
  \int_\Omega |\rmD z|<\infty}$ and we set $\calE(t,\phi,z)=+\infty$ for
$z\not\in \calZ$.

Hence, after minimization with respect to $\phi$ we still have
$\calI(t,z)\geq \gamma +\sigma \int_\Omega |\rmD z|$. This term
provides for $\calR_R^t$ (cf.\ \reff{eq:3sc}) an a priori bound in
$\mathrm{BV}(\Omega,\R^p)$ and hence we conclude compactness in
$X=\rmL^1(\Omega,\R^p)$. Under the usual additional conditions for
the elastic stored-energy densities $W_j$ we obtain for each
$z_0\in \calZ$ a solution $(\phi,z)$ with $\phi \in
\rmL_{\rmw}^\infty(\left]0,T\right[,\rmW^{1,2}(\Omega,\R^d))$ and
$z\in \mathrm{BV}([0,T],\rmL^1(\Omega,\R^p))\cap
\rmL_{\rmw*}^\infty (\left]0,T\right[,  \mathrm{BV}(\Omega,\R^p))
$ with $z(t)\in \calZ$ for all $t\in [0,T]$, see \cite{Main02?}.

\subsection{A delamination problem}
\label{ss:delam}\vskip -5mm \hspace{5mm}

Here we give a simple model for rate-independent delamination and
refer to \cite{KoMiRo02?RIAD} for a better model and the detailed
analysis.

Consider a body $\Omega\subset \R^d$ which is given by an open,
bounded, and path-connected domain. Assume that the interior of
the closure of $\Omega$ differs from $\Omega$ by a finite set of
sufficiently smooth hypersurfaces $\Gamma_j$, $j=1,\ldots,n$. This
means that with $\Gamma := \bigcup^n_{j=1}\Gamma_j$ we have
$\mathrm{int}(\mathrm{cl} (\Omega)) = \Omega \cup \Gamma$.  We
assume that the two sides of the body are glued together along
these surfaces and that the glue is softer than the material
itself. Upon loading, some parts of the glue may break and thus
lose its effectiveness. The remaining fraction of the glue which
is still effective is denoted by the internal state function
$z:\Gamma \to [0,1] $.

We let $\calZ=\set{z:\Gamma \to [0,1]}{ z \text{ measurable}}\subset
X=\rmL^1(\Gamma)$.  The
dissipation distance $\calD(z_0,z_1)$ is proportional to the amount of glue
that is broken from state $z_0$ to state $z_1$:
\[\ts
\calD(z_0,z_1)=c_\calD \int_\Gamma z_0(y){-}z_1(y) \dd a(y)
  \text{ for }z_0\geq z_1 \quad \text{and }+\infty \text{ else}.
\]
Here we explicitly forbid the healing of the glue by setting $\calD$ equal
$\infty$, if $z_0\not\geq z_1$.

The energy is given by the elastic energy in the body, the elastic energy in
the glue, and the potential of the external loadings:
\[\ts
\calE(t,\phi,z)=\int_\Omega W(\rmD\phi)\dd x +\int_\Gamma z(y)
Q(y,\jump{\phi}(y))\dd a(y)-\langle \ell_\ext(t),\phi\rangle,
\]
where for $y\in \Gamma$ the vector $\jump{\phi}(y)$ denotes the jump of the
deformation $\phi$ across the interface $\Gamma$ and $Q(y,\cdot)$ is the
potential for the elastic properties of the glue.

For simplicity we assume further that $W$ provides linearized
elasticity and $Q$ is quadratic is well, then there is a unique
minimizer $\phi=\Phi(t,z)\in \rmH^1(\Omega,\R^d)$ of
$\calE(t,\cdot ,z)$. It can be shown that the mapping
$\Phi(t,\cdot):\calZ\subset \rmL^1(\Gamma) \to
\rmH^1(\Omega,\R^d)$ is compact, which implies that the functional
$\calI(t,\cdot):\calZ\to \left[0,\infty\right[$ is s-weakly
continuous with respect to the $\rmL^1$-topology on $\calZ$. For
the latter argument it is essential that $z$ appears only linearly
in the definition of $\calE(t,\phi,z)$. Theorem \reff{t:3.3} with
condition (a) provides the existence of solutions.

\subsection{Elasto-plasticity}\vskip -5mm \hspace{5mm}

The above theory can be applied to linearized elasto-plasticity,
see \cite{CaHaMi02NCPM,Miel01?FELG}. Here we want to report on
recent results concerning elasto-plasticity with finite strain.
However, for this application the abstract existence theory is not
yet available.

Elasto-plasticity with  finite strains is based on the
multiplicative decomposition of the deformation gradient
$F=\rmD\phi$ in the form $\rmD\phi=F_\mathrm{elast}P^{-1}$ where
the plastic transformation $P$ lies in the Lie group
$\mathrm{SL}(d)=\set{P\in \R^{d\ti d}}{\det P=1}$. The internal
variable has the form $z=(P,p)\in Z$ where $p\in \R^k$ denotes the
hardening parameters. We refer to
\cite{CaHaMi02NCPM,Miel01?FELG,Miel02?EFME} for mechanical
motivations and mathematical details. For simplicity, we mention
here only the case without hardening where $z=P\in \mathrm{SL} (d)
=: Z$ and refer to \cite{HaMiMi02?DDME,Miel02?EFME} for more
general cases.

The important point in finite-strain elasto-plasticity is that the
dissipation distance must be invariant under previous plastic
deformations, i.e., $D(QP_0,QP_1)=D(P_0,P_1)$ for all $Q\in
\mathrm{SL}(d)$. Equivalently, the infinitesimal metric
$\Delta:\rmT Z\to [0,\infty]$ is left-invariant, i.e.,
$\Delta(P,\dot P)= \Delta(I,P^{-1}\dot P)$.  This implies that the
dissipation distance is characterized by a norm $\Delta(I,\cdot)$
on $\mathrm{sl}(d) = \rmT_I \mathrm{SL} (d)$ and that $D(P_0,P_1)$
behaves logarithmically in $P^{-1}_0P_1$ which introduces strong
geometric nonconvexities. So far, even the solution of the
incremental problem (IP) is not understood completely. Even in
simple cases one has to expect non-attainment in (IP), which leads
to the formation of microstructure.  The easiest way to see the
problems is to study the reduced energy density $\Psi^\red$ in
\reff{e4.2}. If this density is not quasi-convex, then there are
loadings such that (IP) has no solution and relaxation techniques
have to be employed, cf.\ \cite{Miel02?RYMM}. \vspace{-1mm}

\paragraph*{Acknowledgments}
This work was partially joint work with Florian Theil, Toma\v{s}
Roub\'{\i}\v{c}ek,  Andreas
Mainik and Michal Ko\v{c}vara. The research was partially supported by DFG
through the SFB 404 \textit{Multifield Problems in Continuum Mechanics}.

\label{lastpage}

\end{document}